\documentclass{agtart_a}
\pdfoutput=1

\usepackage{graphicx,pinlabel,xypic}
\xyoption{all}

\title{Unoriented topological quantum field theory and link homology}
\author{Vladimir Turaev}
\givenname{Vladimir}
\surname{Turaev}
\address{Institut de Recherche Mathematique Avancee\\\newline
7 rue Rene Descartes\\
67000 Strasbourg\\France}
\email{turaev@math.u-strasbg.fr}
\urladdr{}

\author{Paul Turner}
\givenname{Paul}
\surname{Turner}
\address{School of Mathematical and Computer Sciences\\
Heriot-Watt University\\\newline
Edinburgh EH14 4AS\\Scotland}
\email{paul@ma.hw.ac.uk}
\urladdr{}

\volumenumber{6}
\issuenumber{}
\publicationyear{2006}
\papernumber{39}
\startpage{1069}
\endpage{1093}

\doi{}
\MR{}
\Zbl{}


\received{1 September 2005}
\revised{16 January 2006}
\accepted{2 May 2006}
\published{24 August 2006}
\publishedonline{24 August 2006}
\proposed{}
\seconded{}
\corresponding{}
\editor{Martin Scharlemann}
\version{}

\arxivreference{math.GT/0506229}

\let\xysavmatrix\xymatrix
\def\xymatrix{\disablesubscriptcorrection\xysavmatrix}
\AtBeginDocument{\let\bar\wbar\let\tilde\wtilde\let\hat\what}
\AtBeginDocument{\let\colon\colon\thinspace}
\AtBeginDocument{\makeop{Id}}
\AtBeginDocument{\newcommand{\ob}[1]{\wwbar{#1}}}
\AtBeginDocument{\makeop{pt}}

\makeatletter
\def\cnewtheorem#1[#2]#3{\newtheorem{#1}{#3}[section]
\expandafter\let\csname c@#1\endcsname\c@thm}
\makeatother

\newtheorem{thm}{Theorem}[section]
\cnewtheorem{lem}[thm]{Lemma}
\cnewtheorem{cor}[thm]{Corollary}
\cnewtheorem{prop}[thm]{Proposition}
\theoremstyle{definition}
\cnewtheorem{defn}[thm]{Definition}
\cnewtheorem{exe}[thm]{Example}

\makeautorefname{defn}{Definition}
\makeautorefname{exe}{Example}

\newcommand{\bZ}{\mathbb{Z}}
\newcommand{\bQ}{\mathbb{Q}}
\newcommand{\bR}{\mathbb{R}}
\newcommand{\bF}{\mathbb{F}}
\newcommand{\bP}{\mathbb{P}}

\newcommand{\cC}{\mathcal{C}}

\newcommand{\ra}{\rightarrow}
\newcommand{\ot}{\otimes}

\newcommand{\FF}{\bF_2}

\newcommand{\kom}{\mathcal{K}om}
\newcommand{\mat}{\mathcal{M}at}
\newcommand{\modR}{\text{Mod}_R}
\newcommand{\br}[1]{[\mskip-.7mu[ {#1} ]\mskip-.7mu]}
\newcommand{\bri}[2]{[\mskip-.7mu[ {#1} ]\mskip-.7mu]^{{#2}}}
\newcommand{\rptwo}{\bR P^2}

\newcommand{\ucob}[1]{\mathcal{U}\mathcal{C}ob({#1})}
\newcommand{\ucobmodr}[1]{{\ucob {#1}}_{/r}}

\begin{document}

\begin{asciiabstract}
We investigate link homology theories for stable equivalence classes of
link diagrams on orientable surfaces. We apply (1+1)-dimensional
unoriented topological quantum field theories to Bar-Natan's geometric
formalism to define new theories for stable equivalence classes.
\end{asciiabstract}

\begin{htmlabstract}
We investigate link homology theories for stable equivalence classes of
link diagrams on orientable surfaces. We apply (1+1)&ndash;dimensional
unoriented topological quantum field theories to Bar-Natan's geometric
formalism to define new theories for stable equivalence classes.
\end{htmlabstract}

\begin{abstract}
We investigate link homology theories for stable equivalence classes of
link diagrams on orientable surfaces. We apply $(1{+})1$--dimensional
unoriented topological quantum field theories to Bar-Natan's geometric
formalism to define new theories for stable equivalence classes.
\end{abstract}

\maketitle

\section{Introduction}
In this paper we consider link diagrams on orientable surfaces up to the
relation of {\em stable equivalence\/}, that is up to homeomorphisms of
surfaces, Reidemeister moves  and the addition or subtraction of handles
disjoint from
the diagram. Stable equivalence classes of link diagrams have an
equivalent formulation in terms of so called ``virtual'' link diagrams
pioneered by Kauffman (see for example the review articles
Kauffman and Manturov \cite{kauffmanmanturov} and Fenn, Kauffman and Manturov \cite{fennkauffmanmanturov} and references
therein). Many constructions for
link diagrams on $\bR^2$ can be reproduced for stable equivalence
classes of link diagrams on surfaces, for example, one can define the
Jones polynomial.

In his seminal work \cite{khovanov} Khovanov provided a new way to
look at the Jones polynomial of links in $\bR^3$ interpreting it as the Euler
characteristic of a homology theory. This
approach provides an invariant that is not only stronger than the
Jones polynomial but also has nice functorial properties with respect
to link cobordisms. The Jones polynomial may be viewed as a state sum
over the $2^n$ ``smoothings'' of an $n$--crossing link diagram on
$\bR^2$, with each state making a contribution to the polynomial. Each
smoothing is a collection of circles being the result of resolving
each crossing in one of two possible ways (indicated in \fullref{fig:smoothings}). Khovanov's insight was to associate a graded
vector space to each smoothing and arrange these so that 
a certain topological quantum field theory can be used to
define a differential leading to a chain complex associated to the link
diagram. Remarkably the homotopy type of this complex is independent
of the chosen diagram so its homology is an invariant. Bar-Natan has
written a wonderful exposition of all in
\cite{barnatan_khovanov}. The functorial properties of this
construction were conjectured by Khovanov \cite{khovanov}, proven by Jacobsson \cite{jacobsson} and later
reproved in more generality by Bar-Natan \cite{barnatan}.

Khovanov's work does not immediately extend to link diagrams on
arbitrary orientable surfaces. One difficulty arising is the
following. For classical links Khovanov's complex is constructed by
organising the $2^n$ smoothings into a ``cube'' with each edge making
a contribution to the differential. These edges join smoothings with
different numbers of circles and the corresponding term of the
differential may be seen as fusing two circles into one circle or
splitting one circle into two. For links on surfaces this is no longer
the case, and there may be edges joining smoothings with the same
number of circles. This arises when the surface is part of
the structure (as in the work of by Asaeda, Przytycki and Sikora
\cite{asaeda}) and also for stable equivalence classes. In the latter
case Manturov \cite{manturov} solves this problem over the two element
field $\FF$ by setting the differential corresponding to such cube
edges to be zero.  The motivation for the current work was to find
other possibilities for stable equivalence classes.

Bar-Natan's geometric complex \cite{barnatan} which is locally
defined can be extended to diagrams on surfaces, but at the price of admitting
nonorientable cobordisms. In this paper the key new idea is to define the
notion of {\em unoriented\/} topological quantum field theory, which can
then be applied to the geometric picture to obtain a link homology
theory for stable equivalence classes of diagrams on surfaces.

Here is an outline of the paper. In \fullref{sec:utqft} we define
the notion of unoriented topological quantum field theory and study
the underlying algebraic structure which includes a Frobenius algebra
with additional structure. In dimension $1{+}1$ we classify isomorphism
classes of unoriented TQFTs in terms of these Frobenius algebras. We include analysis of the particular case of rank-two theories, which will be used later to give link homologies. In \fullref{sec:complex} we follow Bar-Natan in defining a complex of
cobordisms $\br D$ associated to a link diagram $D$ on a closed
orientable surface. By introducing Bar-Natan's sphere, torus and {\em 4--Tu\/}
relations we obtain an invariant of the link diagram on the surface.
In \fullref{sec:homology} we
apply an unoriented topological quantum field theory to the complex of
cobordisms defined in \fullref{sec:complex} to obtain a complex of modules. Taking homology
yields a link homology. Our main theorem is the following.

\medskip
{\bf \fullref{thm:vlh}}\qua
{\sl Let $R$ be a commutative ring with unit. Let  $\lambda,
\mu,\beta, a, t\in R$ such that $a$ is invertible and the following relations are satisfied:
\begin{eqnarray*}
 &\mu\beta= \lambda\beta = 0,&\\
& 2a\lambda\mu - a^2\mu^2\lambda^2 -a^2\mu^4t = 2.&
\end{eqnarray*}
Then the $5$--tuple $(\lambda, \mu,\beta,a, t)$ defines (the isomorphism
class of) a link homology theory for stable equivalence classes of
diagrams on orientable surfaces.}
\medskip

We end the paper by discussing a few examples.



\section{Unoriented topological quantum field theories}\label{sec:utqft}
In this section we define the
notion of unoriented topological quantum field theory (briefly,
TQFT). In the $1{+}1$--dimensional case there is a classification in terms of
Frobenius algebras with additional structure.

All the constructions can be formulated in the smooth
and PL categories but we prefer to use the language of topological
manifolds. Thus, by a manifold we shall mean a topological manifold. By
a cobordism, we mean a triple $(W,M_0,M_1)$ where $W$ is a compact
manifold whose boundary is a disjoint union of two closed manifolds
$M_0,M_1$ (possibly void). A homeomorphism of cobordisms $(W,M_0,M_1) \ra (W^\prime,M_0^\prime,M_1^\prime)$ is  a homeomorphism $W\ra W^\prime$ sending $M_0$ homeomorphically to $M_0^\prime$ and  $M_1$ homeomorphically to $M_1^\prime$. 

\subsection{The definition}
Fix an integer $d\geq 0$ and a commutative ring with unity $R$.  We
define in this subsection the notion of $(d{+}1)$--dimensional unoriented topological quantum field
theory. Such a TQFT takes values in the category of
projective $R$--modules of finite type (direct summands of $R^n$
with $n=0,1,...$) and $R$--linear homomorphisms. In the case where $R$
is a field, projective $R$--modules of finite type are just
finite-dimensional vector spaces over $R$.

\begin{defn}\label{def:tqft}
A  {\it $(d{+}1)$--dimensional unoriented TQFT $(A,\tau)$\/} assigns a
projective $R$--module  of  finite type $A_M$ to any closed
$d$--dimensional
 manifold  $M$, an
 $R$--isomorphism $f_{\#}\colon A_M\to A_{M'}$ to any  homeomorphism of
$d$--dimensional  manifolds $f\colon M\to M'$, and an $R$--homomorphism
$\tau(W)\colon  A_{M_0}\to A_{M_1} $ to any $(d{+}1)$--dimensional
cobordism     $(W,M_0,M_1)$. These modules and homomorphisms
should satisfy the following  seven axioms.

(1)\qua For any homeomorphisms of closed  $d$--dimensional manifolds
$f\colon M\to M'$ and $f'\colon M'\to M''$, we have $(f'f)_{\#}=f'_{\#}f_{\#}$. The
isomorphism $f_{\#}\colon A_M\to A_{M'}$ is invariant under   isotopies
of  $f$.

(2)\qua For any disjoint closed $d$--dimensional manifolds $M,N$, there is
an isomorphism $A_{M\amalg N}=A_M\otimes A_N$, natural with respect to
homeomorphisms, where $\otimes=\otimes_R $ is the tensor product over
$R$.

(3)\qua  $  A_{\emptyset}=R$ where the empty set is considered as a
$d$--dimensional manifold (for any $d$).

(4)\qua The  homomorphism  $\tau $ associated with cobordisms is
natural with respect to homeomorphisms of cobordisms.

(5)\qua If a $(d{+}1)$--dimensional cobordism $W$ is a disjoint union of
cobordisms $W_1,W_2$, then under the identifications in axiom (2),
$\tau(W)=\tau(W_1) \otimes
\tau(W_2)$.

(6)\qua For a cobordism     $(W,M_0,M_1)$ obtained from 
$(d{+}1)$--dimensional cobordisms $(W_0,M_0,N)$ and $(W_1,N',M_1)$ by
gluing along a  homeomorphism $f\colon N\to N'$, 
$$\tau(W)= \tau(W_1)\circ f_{\#} \circ \tau(W_0)\colon A_{M_0}\to A_{M_1}.$$

(7)\qua For any closed  $d$--dimensional manifold
 $M$,
we have   
$$\tau(M\times [0,1], M\times 0, M\times1)=\Id\colon  {A_M}\to {A_M}
$$ 
where  we identify $M\times 0$ and
$M\times 1$ with $M$ in the obvious way.
\end{defn}

The pair $(A,\tau)$ is {\em isomorphic\/} to another
$(A^\prime,\tau^\prime)$ if for each closed $d$--manifold $M$ there is
an isomorphism $\eta_M\colon A_M\ra A^\prime_M$, natural with respect to
homeomorphisms and cobordisms, multiplicative with respect to disjoint
union and such that $\eta_\emptyset = \Id_R$.

This definition is essentially the one of Atiyah \cite{atiyah} but all
references to orientations of manifolds are suppressed. The axioms
(1)--(7) constitute a special case of a detailed axiomatic definition
of TQFT's given in Turaev \cite{turaev}, Chapter III in a framework of
so-called space-structures.  For more on the naturality in axioms 2
and 4, the reader is referred to \cite[p.\ 121]{turaev}; otherwise a
knowledge of \cite{turaev} will not be required and the language of
space-structures will not be used.

In dimension $1{+}1$ a related idea has already appeared in Alexeevski and Natanzon
\cite{alexeevskinatanzon} in the context of open--closed field
theory. The definition presented there of a Klein topological field
theory is, however, rather different from that above, justifying a
separate treatment here.

Note that the tensor product in axiom (2) is {\it unordered\/} so
that the modules $ A_M\otimes A_N$ and $ A_N\otimes A_M$ are the same
and not merely isomorphic. The unordered tensor product of a finite
family of modules is obtained by considering all possible orderings of
this family, forming the corresponding ordered tensor products and
then identifying the resulting modules along the obvious isomorphisms
induced by permutations. We have to use here the unordered tensor
product since we want to apply axiom (2) to nonconnected $d$--manifolds
without fixing any order in the set of connected components.

\subsection{Extended Frobenius algebras}
 We recall that a {\em commutative Frobenius algebra over $R$\/} is a
 unital, commutative $R$--algebra $V$ which as an $R$--module is
 projective of finite type, together with a module
 homomorphism $\epsilon\colon V \ra R$ such that the bilinear form
 $\langle - , - \rangle\colon V\otimes V \ra R$ defined by $\langle v,
 w \rangle = \epsilon (vw)$ is nondegenerate, ie, the
 adjoint homomorphism $V\ra V^*$ is an isomorphism. The homomorphism
 $\epsilon$ is called the co-unit and it is useful to define a
 coproduct $\Delta\colon V \ra V\ot V$ by $\Delta (v) = \sum_i
 v_i^\prime \ot v_i^{\prime\prime}$ being the unique element such that
 for all $w\in V$, $vw = \sum_i v_i^\prime \langle
 v_i^{\prime\prime},w\rangle$.

The map $\Delta\colon V \ra
V\otimes V$ satisfies 
\[
(\Id \otimes m ) \circ (\Delta \otimes \Id) = \Delta \circ m = (m\otimes \Id) \circ (\Id \otimes \Delta)
\]
\[
(\Id \otimes \epsilon)\circ \Delta = \Id = (\epsilon \otimes \Id) \circ \Delta
\leqno{\hbox{and}}
\]
where $m\colon V\otimes V \ra V$ is the multiplication in $V$ and
$i\colon R \ra V$ is the unit map.

An unoriented $(d{+}1)$--dimensional TQFT $(A,\tau)$ has an underlying
Frobenius algebra equipped with an involution of Frobenius algebras as we now discuss.

For the underlying $R$--module we take $V = A_{S^d}$, where $S^d$ is
the standard unit sphere in Euclidean space $\bR^{d+1}$. Note that
any other sphere $S$ in $\bR^{d+1}$ can be related to
the standard sphere $S^d$ by parallel translations and homotheties
which allows us to canonically identify the corresponding vector
spaces (courtesy of axiom (1) in \fullref{def:tqft}).

The unit $i$ and the co-unit $\epsilon$ come from the
$(d{+}1)$--dimensional unit ball $B$ in $\bR^{d+1}$ viewed as a cobordism
from the empty set to $S^d$ and as a cobordism from $S^d$ to the empty
set respectively.

The product is obtained by considering a ball $W_{00}$ of radius three
in $\bR^{d+1}$ with two interior balls of radius 1 removed. This can
be viewed as a cobordism from the two internal spheres to the external
sphere which, using the canonical identifications above, gives a map
$m\colon V\otimes V \ra V$. Although not necessary as part of the
definition we note that the coproduct has geometric interpretation by
viewing $W_{00}$ as a cobordism from the sphere of radius three to the
two internal spheres.

The involution $\phi\colon V \ra V$ is induced by an orientation-reversing
homeomorphism $\chi\colon S^d\ra S^d$. They are all isotopic so $\phi$ is well
defined.

\begin{prop}
The $R$--module $V$ equipped with structure maps $m,i, \epsilon$
defined above is a Frobenius algebra and $\phi$ is an involution of
Frobenius algebras.
\end{prop}

\begin{proof}
It follows from the familiar oriented case that $(V,m,i,\epsilon)$ is a Frobenius algebra.

We observe that $\phi^2=\Id$ since $\chi\circ \chi$ is isotopic to the identity on $S^d$.

Next we claim that $\phi$ is a homomorphism of Frobenius
algebras. Note that an orienta\-tion-reversing homeomorphism $S^d\ra
S^d$ extends to $B$, hence by axiom (4) in \fullref{def:tqft},
$\phi$ preserves the unit and counit. There is also an
orientation-reversing homeomorphism $W_{00} \ra W_{00}$ mapping each
boundary component to itself showing that $\phi \circ m \circ (\phi\ot
\phi) = m$, and hence since $\phi^2=\Id$, we see $m\circ (\phi\ot \phi) =
\phi \circ m$.\proved
\end{proof}

One further piece of structure can be identified in the form of an
element $\theta \in V$. Consider the punctured projective
space $\bP$ of dimension $d+1$ viewed as a cobordism between the
empty set and $\partial \bP$. Now we identify $\partial \bP$ with $S^d$
via a map $f\colon \partial \bP \ra S^d$ and set $\theta = f_\#\circ
\tau(\bP)(1)$. There is an involution $T$ on $\bP$ given by the
negation of one coordinate which by axiom (4) in \fullref{def:tqft} gives $T_\# ( \tau(\bP)(1)) = \tau(\bP)(1)$. 
Now $T$ reverses the orientation of the boundary and thus if
$f^\prime$ is another homeomorphism $\partial \bP \ra S^d$
nonisotopic to $f$ then $f^\prime \circ T$ is isotopic to $f$ and we
have
\[
f^\prime_\# ( \tau(\bP)(1)) = f^\prime_\#(T_\#( \tau(\bP)(1))) = f_\#(
\tau(\bP)(1)).
\]
This shows that $\theta$ is well defined.

\begin{lem}
$\phi(\theta) = \theta$.
\end{lem}
\begin{proof}
 The map $T\colon \bP \ra \bP$ defined above reverses the
orientation on $\partial \bP$ so by axiom (4) in \fullref{def:tqft} we have $\phi(\theta) = \theta$.
\end{proof}

When $d$ is odd there is a stronger relation between $\phi$ and $\theta$.

\begin{prop}\label{prop:ptv}
If $d$ is odd then $\phi(\theta v) = \theta v$ for all $v\in V$.
\end{prop}

\begin{proof}
Consider two disjoint balls $B_1$ and $B_2$ in $\bP^{d+1}$ in the
complement of $\bP^d$ in $\bP^{d+1}$. Let $S_i = \partial B_i$,
$i=1,2$ and consider the cobordism $(W,S_1,S_2)$ where $W= \bP^{d+1}
- \text{Int}(B_1\sqcup B_2)$. It is clear that $W$ is a connected sum
$C \,\#\, \bP^{d+1}$ where $C$ is homeomorphic to $S^d\times I$. As the
$d$--sphere separating $W$ into two connected summands we take the
boundary of a regular neighbourhood of $\bP^d$ in $\bP^{d+1}$. Now
choose an arbitrary orientation of $S_1$ and endow $S_2$ with the
orientation such that $\tau(C)\colon V\cong A_{S_1} \ra A_{S_2}\cong
V$ is the identity homomorphism. We can then use the decomposition of
$W$ to compute $\tau(W)(v) = \theta v$.

Since $d$ is odd, and hence $\bP^{d+1}$ unorientable, there is a
homeomorphism $T\colon W \ra W$ which reverses the orientation of
$S_2$ while leaving the orientation of $S_1$ unchanged. This is
obtained by moving $B_2$ by an isotopy in $\bP^{d+1}$ back to itself
passing through $\bP^d$ exactly once. Since $T$ reverses the
orientation of $S_2$, using the identifications above we have
$T|_{S_2} = \phi\colon V \ra V$. On $S_1$ we have $T|_{S_1} = \Id\colon
V \ra V$. By axiom (4) in \fullref{def:tqft} the following
diagram commutes.
\[
\xymatrix{
V\ar[d]_{\Id}  \cong A_{S_1} \ar[r]^{\tau(W)}& A_{S_2}\cong  V \ar[d]^\phi \\
V  \cong A_{S_1} \ar[r]^{\tau(W)}& A_{S_2}\cong  V 
}
\]
Chasing around this diagram gives $\phi(\theta v) = \theta v$ as required.
\end{proof}

In fact in dimension $1{+}1$ the assignment of the structure above
induces a bijection from isomorphism classes of TQFTs to a special
class of Frobenius algebras.

\begin{defn}\label{def:efa}
An {\em extended Frobenius algebra\/} is a Frobenius algebra  $(V,m,i,\epsilon)$ together with an involution of Frobenius algebras $\phi\colon V \ra V$ and an element $\theta\in V$ satisfying the following two axioms.
\begin{enumerate}
\item $\phi(\theta v)=\theta v$, for all $v\in V$.
\item  $m(\phi\otimes \Id)(\Delta(1)) = \theta^2$.
\end{enumerate}
\end{defn}

Two extended Frobenius algebras  $(V,\theta,\phi)$ and $(V^\prime,
\theta^\prime,\phi^\prime)$ are {\em isomorphic\/} if there exists an
isomorphism of Frobenius algebras $g\colon V\ra V^\prime$ such that
$g(\theta) = \theta^\prime$ and $g\circ \phi = \phi^\prime \circ g$.

We note that in \cite{alexeevskinatanzon} such algebras appear as
``structure algebras'' in which the part coming from open boundaries
is trivial.


\begin{exe}
Let $N$ be odd and let $V=\bR[x]/x^N$. Taking $\phi=\Id$ and $\theta 
= \sqrt{N} x^{\frac{N-1}{2}}$ we see that $(V, \phi,\theta)$ is an extended Frobenius algebra. In this case $\phi$ is obviously an involution of Frobenius algebras satisfying axiom (1). Axiom (2) is reduced to $m(\Delta(1)) = \theta^2$ which holds for our choice of $\theta$ since $m(\Delta(1)) = Nx^{N-1}$.
\end{exe}

\begin{exe}
Let $V$ be the two dimensional vector space over $\bF_2$ on generators $1$ and $x$. This may be given the structure of  Frobenius algebra by defining $i(1) = 1$, $x^2 = x$, $\epsilon(1) = 0$, $\epsilon (x) =1$. Taking  $\theta = 0$ and $\phi(1) = 1$, then $\phi(x) = 1+x$ turns this into an extended Frobenius algebra.
\end{exe}

There are some useful elementary consequences of the definition.

\begin{lem}\label{lem:useful}
{\rm(i)}\qua $m(\phi\otimes \Id)(\Delta(v)) = m(\phi\otimes \Id)(\Delta(1))v$ for all
$v\in V$.

{\rm(ii)}\qua  $m(\phi\otimes \Id)(\Delta(v)) = \theta^2v$ for all $v\in V$.

{\rm(iii)}\qua $m(\Delta(\theta)) = \theta^3$.
\end{lem}

\begin{proof}
To prove (i) it suffices to notice that $\Delta(v) = (v\ot 1)\Delta (1)$.
Part (ii) follows immediately from (i).
For (iii) we have
\[
m(\Delta(\theta)) = m((\theta\ot 1)\Delta(1)) = m(\phi\ot 1)((\theta\ot 1)\Delta(1)) = m(\phi\ot 1)(\Delta(\theta)) = \theta^3.
\proved \]  
\end{proof}

\begin{prop}\label{prop:tqftfrob}
Isomorphism classes of unoriented $1{+}1$--dimensional TQFTs over $R$ are
in bijective correspondence with isomorphism classes of extended Frobenius
algebras over $R$.
\end{prop}

\begin{proof}
We have already shown that an unoriented $1{+}1$--dimensional TQFT has an
underlying Frobenius algebra over $R$ and comes equipped with an
involution of Frobenius algebras $\phi$. Furthermore we have defined
the element $\theta$ and shown in \fullref{prop:ptv} that
$\phi(\theta v) = \theta v$ for all $v\in V$. It remains to show axiom
(2) in the definition of extended Frobenius algebra. To see this we
glue the ends of an oriented punctured cylinder together in a
nonorientation-preserving way, regarding the result $W$ as a
cobordism $(W,\emptyset, \partial W)$. We can use the original
orientation of the cylinder to orient $\partial W$ and hence identify
$A_{\partial W}$ with $V$. One can easily compute $\tau(W)\colon R \ra
V$ to be the map taking 1 to $m(\phi \ot \Id)(\Delta(1))$. On the other
hand $W$ is a punctured Klein bottle and hence homeomorphic to the
punctured connected sum of two projective planes. Thus $\tau(W)(1)
=\theta^2$.  It is clear that the isomorphism class of this algebra
depends only on the isomorphism class of the TQFT.

We now show that every extended Frobenius algebra $(V,\theta,\phi)$
gives rise to an unoriented TQFT. Given a connected closed $1$--manifold
$M$ we define the set $A_M$ by
\[
A_M = \{(\gamma,v) \mid \gamma\colon S^1 \ra M \text{ a homeomorphism}, v\in V\}/\sim
\]
\begin{align}
(\gamma,v) \sim (\gamma^\prime,v^\prime) \text{ iff }  & \text{{\it either \/}}\gamma \text{ is isotopic to } \gamma^\prime \text{ and } v=v^\prime 
\tag*{\hbox{where}} \\
& \text{{\it or \/}}\gamma \text{ is not isotopic to } \gamma^\prime \text{ and } v=\phi(v^\prime).\notag
\end{align}
Note that there are two isotopy classes of homeomorphisms $S^1\ra M$. Now pick any homeomorphism $h\colon S^1\ra M$ and define a map (of sets) $\tilde h\colon A_M\ra V$ by
\[
\tilde h(\gamma,v) = \begin{cases}
v & \text{ if $\gamma$ is isotopic to $h$}\\
\phi(v) & \text{ else.}
         \end{cases}
\]
This map is a bijection and moreover, since $V$ is an $R$--module, we
can use it to turn $A_M$ into an $R$--module. It is easy to check that
the $R$--module structure is independent of the isotopy class of
$h$. Thus to a closed $1$--manifold $M$ we have assigned an $R$--module
$A_M$. If $M$ is not connected we simply write $M=M_1\sqcup \cdots
\sqcup M_k$ and let $A_M = \bigotimes A_{M_i}$. We define $A_\emptyset = R$.

Now let $f\colon M\ra M^\prime$ be a homeomorphism of connected
$1$--manifolds (up to isotopy there are two choices). Define $f_\#\colon
A_M \ra A_{M^\prime}$ by $f_\#(\gamma,v) = (f\circ \gamma,v)$. One can
check this is a well defined isomorphism of $R$--modules. This extends to homeomorphisms of 
nonconnected manifolds by multiplicativity.

We now wish to define $\tau(W)\colon A_{M_0} \ra A_{M_1}$ for a
cobordism $(W,M_0,M_1)$. We separate the cases of
orientable and nonorientable surfaces. 

Suppose first that $W$ is orientable and connected. Choose an
orientation for $W$ and decompose it into basic pieces (caps, cups and
pairs of pants).
We can appeal to the usual case of oriented TQFT to define $\tau(W)$
which is independent of the decomposition.  We claim
that $\tau(W)$ is independent of the orientation chosen for $W$. Indeed by
choosing the opposite orientation, each piece in the decomposition
also has the opposite orientation. Consider, for example a pair of
pants $(P,M,N)$ occurring in the decomposition and consider the
following diagram where the top route corresponds to one orientation
and the bottom route to the other.
\[
\xymatrix{
& V^{\ot 2}  \ar[r]^m \ar[dd]_{\phi^{\ot 2}} & V \ar[dr]^{\cong} \ar[dd]^\phi & \\
A_M \ar[ru]^\cong \ar[rd]_\cong & & & A_N\\
& V^{\ot 2}  \ar[r]^m & V \ar[ur]_{\cong} & 
}
\]
Note that if a connected $1$--manifold $M$ is oriented then there is a
canonical identification of $A_M$ with $V$, using the map $\tilde
h\colon A_M \ra V$ defined by an orientation-preserving homeomorphism
$h\colon S^1\ra M$ (here $S^1$ is given the anticlockwise
orientation). By reversing the orientation of $M$ the identification
is given by $\phi\circ \tilde h$.  Thus the left and right triangles in
the diagram above commute. Since $\phi$ is a map of Frobenius algebras
the middle square also commutes, hence the two routes give the same
map. Similar arguments hold for the other basic surfaces, showing that
$\tau(W)$ is independent of the orientation. Furthermore, the
properties of oriented TQFTs guarantee that $\tau(W)$ is natural with
respect to homeomorphisms. For nonconnected $W$ we extend the above
multiplicatively.

Suppose now that $W$ is nonorientable and connected. We may present
$W$ as a connected sum of an orientable surface $W^{\text{or}}$ and
$n$ projective planes, $W=W^{\text{or}} \,\#\, n \rptwo$. Note that
$\partial W^{\text{or}}= \partial W$ and that the homomorphism
$\tau(W^{\text{or}})\colon A_{M_0}
\ra A_{M_1}$ is defined by the orientable case above. 
Now choose an identification $A_{M_1}\cong V^{\ot k}$ and define 
$\psi_n\colon A_{M_1}\cong V^{\ot k} \ra V^{\ot k} \cong A_{M_1}$ to be 
the identity on all factors except one where it is multiplication by
$\theta^n$. We define $\tau(W) = \psi_n \circ \tau(W^{\text{or}})$.

A priori this depends on the identification $A_{M_1}\cong V^{\ot k}$,
the choice of factor in $\psi_n$ and the decomposition $W=
W^{\text{or}} \,\#\, n \rptwo$. It follows from the properties of $V$ that
the factor in $\psi_n$ does not matter. Moreover, since $\phi(\theta) = \theta$ the definition is independent of the identification $A_{M_1}\cong V^{\ot k}$. Now suppose that we decompose
$W$ differently as $W=\ob{W}^{\text{or}}
\, \#\, n \rptwo$ then there is a homeomorphism $g\colon \ob{W}^{\text{or}} \ra W^{\text{or}}$ taking each boundary component to itself. On each such
boundary component $g$ is either isotopic to the identity or is an
orientation-reversing homeomorphism. By the naturality of $\tau$ for
orientable surfaces we have the following commutative diagram.
\[
\xymatrix{
A_{M_0} \ar[r]^{\tau(W^{\text{or}})} \ar[d]_{(g|_{M_0})_\#}   & A_{M_1} \ar[d]^{(g|_{M_1})_\#}\\
A_{M_0} \ar[r]^{\tau(\ob{W}^{\text{or}})}    & A_{M_1}
}
\]
We also have the following commutative diagram.
\[
\xymatrix
{
A_{M_1} \ar[r]^{\psi_n} \ar[d]_{(g|_{M_1})_\#}   & A_{M_1} \ar[d]^{(g|_{M_1})_\#}\\
A_{M_1} \ar[r]^{\psi_n}    & A_{M_1}
}
\]
Commutativity is immediate for a component of $M_1$ on which
$g|_{M_1}$ is isotopic to the identity and follows from the
relation $\phi(\theta v) = \theta v$ on a component where
$g|_{M_1}$ is an orientation-reversing homeomorphism. Combining
these two diagrams shows that $\tau(W)$ is independent of the
decomposition above. We may also choose a different number of
projective planes in the decomposition, replacing any three by a torus
and a projective plane. However, in this case we may write
$W=W^{\text{or}}\,\#\, T^2\,\#\, (n-2)\rptwo$ and $\tau(W)$ can be computed as
$\tau(W^{\text{or}})$ multiplied on one factor by
$m(\Delta(1))\theta^{n-2}$. However, by \fullref{lem:useful} we have
$m(\Delta(1))\theta^{n-2} = m(\Delta(\theta))\theta^{n-3} = \theta^n$
showing that $\tau(W)$ is again independent of the decomposition. For a nonconnected cobordism we extend multiplicatively.

What remains is to show that the structure defined above verifies
axioms (1)--(7) in \fullref{def:tqft}. These are all easy to
show with the exception of the gluing axiom (6). We leave the others
as an exercise and focus on proving this axiom, for which the following definition 
is useful.

\begin{defn}
A cobordism $(W_0,M_0,N)$ is {\em nice\/} if for all $(W_1,N^\prime,M_1)$ and homeomorphisms $f\colon N \ra N^\prime$ we have
\[
\tau(W) = \tau(W_1)\circ f_\# \circ \tau(W_0)
\]
where $(W,M_0, M_1)$ is the result of gluing $W_0$ to $W_1$ along $f$.
\end{defn}

It is easy to check the following lemma.

\begin{lem}\label{lem:nice}
(i)\qua If the cobordism $(W_0,M_0,N_0)$ is obtained from $(W_0^\prime,M_0^\prime,N_0^\prime)$ and $(W_0^{\prime\prime},M_0^{\prime\prime} ,N_0^{\prime\prime} )$ by
gluing along a homeomorphism $g\colon N_0^\prime \ra M_0^{\prime\prime}$, and $W_0^\prime$ and
$W_0^{\prime\prime}$ are nice, then $W_0$ is also nice.

(ii)\qua If $(W_0,M_0,N_0)$ and $(W_1,M_1,N_1)$ are oriented and $f\colon N_0 \ra M_1$ is orientation preserving then 
$\tau(W) = \tau(W_1)\circ f_\# \circ \tau(W_0)$.

\end{lem}

Any surface may be decomposed into a composition of cobordisms of the
six basic types (up to ordering) indicated in \fullref{fig:basic}
where the input boundary is always at the top (the picture on the
right at the bottom represents a twice punctured projective plane).

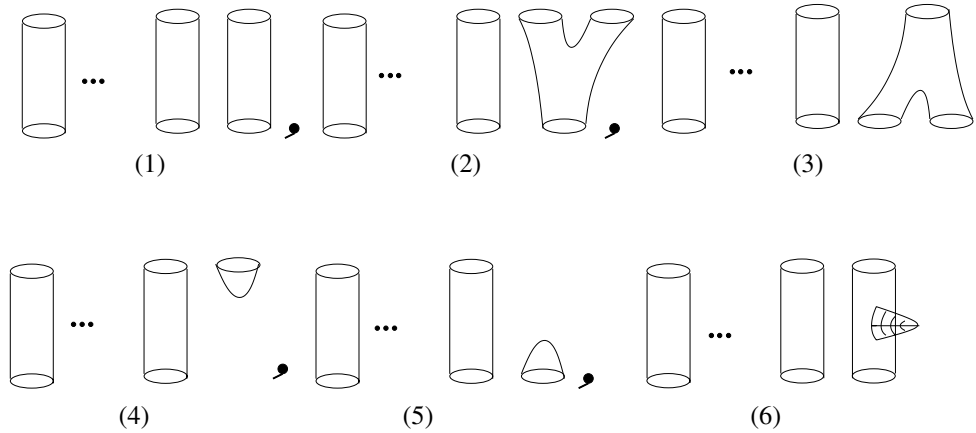
\begin{figure}[ht!]
\labellist
\small\hair 2pt
\pinlabel $(1)$ [t] at 53 109
\pinlabel $(2)$ [t] at 172 109
\pinlabel $(3)$ [t] at 301 109
\pinlabel $(4)$ [t] at 47 14
\pinlabel $(5)$ [t] at 154 14
\pinlabel $(6)$ [t] at 285 14
\endlabellist
\centerline{\includegraphics{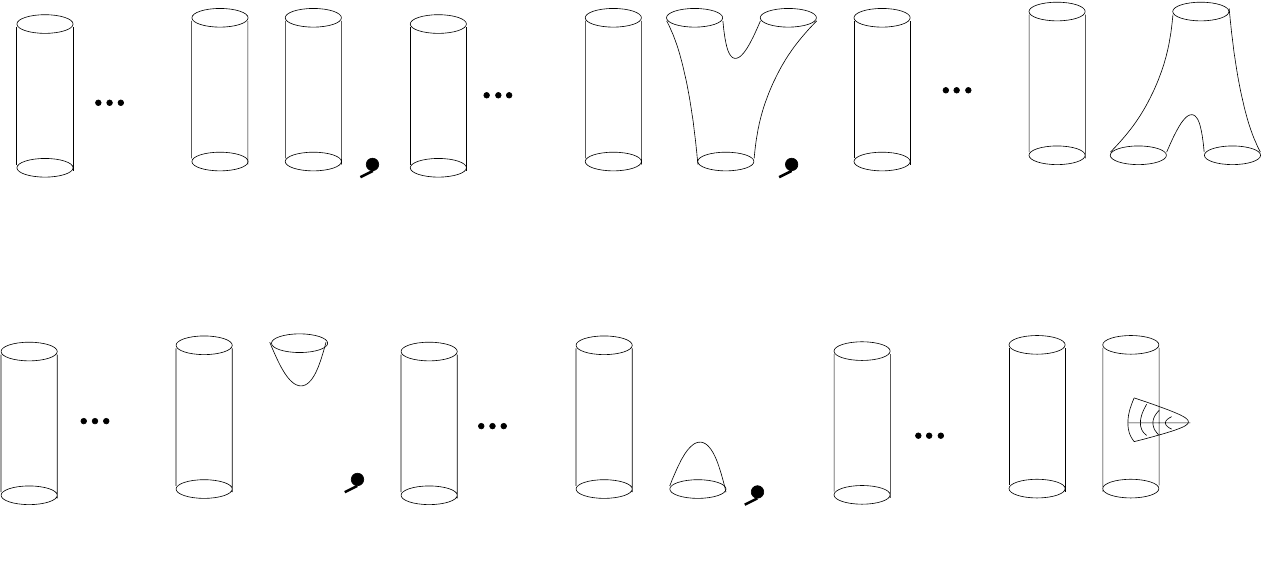}}
\caption{Basic cobordism types}
\label{fig:basic} 
\end{figure}

In order to prove that the gluing axiom holds we will show that any
cobordism of the six types in \fullref{fig:basic} is nice. Thus by \fullref{lem:nice} (i) we
will have the desired result.

First we claim that any cobordism $W_0$ of type (2) is nice. To see
this let $(W_1,N^\prime,M_1)$ be a cobordism which we glue to $W_0$
along a homeomorphism $f\colon N\ra N^\prime$. We can present $W_1$ as
$W_1=W_1^{\text{or}} \,\#\, n\rptwo$ ($n\geq 0$). We now orient
$W_1^{\text{or}}$ and choose orientations for each of the cylinders
and the pair of pants in $W_0$ such that $f$ is orientation
preserving. Let $W_2$ be the result of gluing $W_1^{\text{or}}$ to
$W_0$ along $f$ so that $W=W_2 \,\#\, n\rptwo$. By \fullref{lem:nice}
(ii) we have $\tau(W_2) = \tau(W_1^{\text{or}})\circ f_\# \circ
\tau(W_0)$.  Thus
\begin{align*}
\tau(W) & = \psi_n  \circ \tau(W_2) =  \psi_n \circ (\tau(W_1^{\text{or}})\circ f_\# \circ \tau(W_0)) \\
& =  (\psi_n \circ \tau(W_1^{\text{or}})) \circ f_\# \circ \tau(W_0) = 
 \tau(W_1)\circ f_\# \circ \tau(W_0).
\end{align*}
Similar arguments show that cobordisms of type (1), (4) and (5) are nice.

For a cobordism of type (6) we decompose $W_0$ and
$W_1$ as $W_0 = W_0^{\text{or}} \,\#\, \rptwo$ and $W_1 = W_1^{\text{or}}
\,\#\, n\rptwo$. Now let $W_2$ be the result of gluing $W_1^{\text{or}}$
to $W_0^{\text{or}}$ along $f$. We can write $W = W_2 \,\#\, (n+1)\rptwo$. Thus
\begin{align*}
\tau(W) & = \psi_{n+1}  \circ \tau(W_2) = \psi_{n+1} \circ (\tau(W_1^{\text{or}}) \circ f_\# \circ \tau(W_0^{\text{or}}))\\
 & = (\psi_1 \circ \tau(W_1^{\text{or}})) \circ f_\# \circ (\psi_n\circ \tau(W_0^{\text{or}})) = \tau(W_1) \circ f_\# \circ \tau(W_0).
\end{align*}

The most difficult case is for cobordisms of type (3). We decompose
$W_1$ into the connect sum $W_1^{\text{or}}\,\#\, n\rptwo$ and we orient
$W_1^{\text{or}}$ in an arbitrary way. If there is an orientation of
$W_0$ such that $f$ is an orientation-preserving homeomorphism then we
are done by the same arguments as above. If there is no such
orientation we decompose $W_1$ into the composition of $W_2$ and $W_3$
glued along a (canonical) map $g$ where $W_2$ is of type (2). Since by
the above $W_2$ is nice we have $\tau(W_1) = \tau(W_3)\circ g_\#
\circ \tau(W_2)$. Now let $W_4$ be the result of gluing $W_2$ to $W_0$
along $f$, which results in a number of cylinders and a twice
punctured Klein bottle. We can compute $\tau(W_4) = \psi_2 \circ
\tau(W_4^{\text{or}})$. On the other hand by using \fullref{lem:useful} (ii) we have $\tau(W_4) = \tau(W_2)
\circ f_\# \circ \tau(W_0)$. Finally we note that $W$ is the result of gluing $W_3$ to $W_4$ along $g$ and so 
\begin{align*}
\tau(W) & = \psi_{n+2} \circ   \tau(W^{\text{or}}) 
= \psi_{n+2} \circ (\tau(W_3^{\text{or}}) \circ g_\# \circ
\tau(W_4^{\text{or}}))\\ & = (\psi_{n} \circ \tau(W_3^{\text{or}}))
\circ g_\# \circ (\psi_2 \circ \tau(W_4^{\text{or}})) = \tau(W_3)
\circ g_\# \circ \tau(W_4) \\ & = \tau(W_3) \circ g_\# \circ \tau(W_2)
\circ f_\# \circ \tau(W_0) 
= \tau(W_1) \circ  f_\# \circ \tau(W_0). 
\end{align*}
The gluing axiom now follows since all the pieces in the decomposition are nice.

Thus we have defined an unoriented TQFT. Standard arguments now show
that the isomorphism class of this TQFT depends only on the
isomorphism class of the extended Frobenius algebra. Moreover, the
underlying extended Frobenius algebra is clearly $(V,\theta,\phi)$ so
the construction provides the required inverse.
\end{proof}

We remark that in order to use the underlying Frobenius algebra to
make computations we must choose for each closed $1$--manifold $\Gamma$ an
identification $\gamma \colon S^1 \sqcup \cdots \sqcup S^1 \ra \Gamma$
which gives an isomorphism $\tau(\gamma) \colon V^{\ot r} \ra
A_{\Gamma}$ where $r$ is the number of components of $\Gamma$.

\subsection{Rank-two aspherical theories}
In this subsection we wish to study rank-two unoriented TQFTs
satisfying the condition $\tau(S^2) =0$. For convenience we refer to these as rank-two {\it aspherical\/} unoriented TQFTs.

\begin{lem}
Let $(A,\tau)$ be a rank-two aspherical unoriented TQFT and let $S$ be any sphere and $T$ any torus. Then $\tau(S) = 0$ and $\tau(T) = 2$.
\end{lem}
\begin{proof}
We have a standard torus constructed from discs and the surface
$W_{00}$ (a pair of pants surface) which evaluates to 2 because the
theory is rank-two. It follows from this, the fact that the TQFT is
aspherical and axiom (4) in \fullref{def:tqft} that any sphere evaluates to 0 and any torus evaluates to 2.
\end{proof}

The following is an immediate corollary to \fullref{prop:tqftfrob}.

\begin{prop}
The isomorphism classes of rank-two aspherical unoriented $1{+}1$--dimensional
TQFTs over $R$ are in bijective correspondence with the isomorphism
classes of rank-two extended Frobenius algebras over $R$ satisfying
$\epsilon(i(1)) = 0$, where $\epsilon$ and $i$ are the counit and unit
of the Frobenius algebra.
\end{prop}

We now wish to classify the Frobenius algebras appearing in the
previous proposition. Let $K=\bZ [a,f,t,\lambda,\mu,\beta]/I $ where
$I$ is the ideal generated by
\[
af=1, \quad \mu\beta= \lambda\beta = 0 \quad \text{and} \quad 2a\lambda\mu - a^2\mu^2\lambda^2 -a^2\mu^4t = 2.
\]
Now set $U=K\{1,x\}$ and define multiplication by
\[
11=1, \quad 1x=x1=x \quad \text{and} \quad xx= ( \beta - a \lambda^2 -a \mu^2t)x+t1
\]
and comultiplication by
\begin{gather*}
\Delta (1) = f(1\otimes x + x \otimes 1)  -( \beta - a \lambda^2 -a \mu^2t)f1 \otimes 1,\\
 \Delta
(x) = fx\otimes x + ft1\otimes 1.
\end{gather*}
Define a unit and counit by
\[
\epsilon(1)=0, \quad \epsilon (x) = a \quad \text{and} \quad i(1) = 1.
\]
\[
\theta = \lambda 1 + \mu x \leqno{\hbox{Let}} 
\]
and define $\phi\colon U\ra U$ by 
\[
\phi(1) = 1 \quad \text{and} \quad \phi(x) = \beta 1 + x.
\]

\begin{prop}\label{prop:univefa}
The triple $(U,\theta,\phi)$ defined above is an extended Frobenius
algebra over $K$ satisfying $\epsilon(i(1)) =0$.
\end{prop}

\begin{proof}
The proof is purely computational and hence omitted. The one thing
worth pointing out is that the relations in $K$ imply $2\beta =0$.
\end{proof}

Given another ring with unity $R$ and a ring homomorphism $\psi\colon
K \ra R$ we can use $\psi$ to view $R$ as an $K$--module.  We set
$U_\psi = U\otimes_K R$ and define a comultiplication and counit by
$\Delta \otimes 1$ and $\epsilon \otimes 1$. Since $U$ is a Frobenius
algebra over $K$, it follows that $U_\psi$ is a Frobenius algebra over
$R$. The resulting triple $(U_\psi,
\theta\otimes 1, \phi\otimes \Id)$ is an extended Frobenius algebra
over $R$ and it too satisfies $\epsilon (i(1)) =0$.

The example in \fullref{prop:univefa}
is universal for rank-two theories in the sense of the following proposition.

\begin{prop}\label{prop:univ}
Let $R$ be a commutative ring with unit and let
$(V,\theta^\prime,\phi^\prime)$ be a rank-two extended Frobenius
algebra over $R$ satisfying $\epsilon(i(1)) =0$. Then there exists a ring
homomorphism $K\ra R$ such that $(U_\psi, \theta\ot 1,\phi\ot \Id)$ is
isomorphic to $(V,\theta^\prime,\phi^\prime)$.
\end{prop}

\begin{proof}
Any rank-two Frobenius algebra $V$over $R$ satisfying $\epsilon(i(1))
=0$ is of the following form (see Khovanov \cite{khovanov_frob}). As an
$R$--module we can write $V=R\{1,x\}$ and there are elements $a,f,h,t\in R$
such that $af=1$. In terms of these elements multiplication is defined
by
\[
11=1, \quad 1x=x1=x \quad \text{and} \quad xx= hx+t1
\]
and comultiplication is defined by
\[
\Delta (1) = f(1\otimes x + x \otimes 1)  -hf1 \otimes 1 \quad \text{and} \quad \Delta
(x) = fx\otimes x + ft1\otimes 1.
\]
The unit and counit are given by
\[
\epsilon(1)=0, \quad \epsilon (x) = a \quad \text{and} \quad i(1) = 1.
\]
Suppose now $\theta\in V$ and $\phi\colon V\ra V$ such that
$(V,\theta,\phi)$ is an extended Frobenius algebra. Since $\phi$ is a
map of Frobenius algebras we must have $\phi(1) = 1$.  Now write
\[
\theta = \lambda 1 + \mu x \quad \text{and} \quad 
\phi(x) = \beta 1 + \gamma x.
\] 
We have $\epsilon(x) = \epsilon(\phi(x))$ from which we have
\[
a= \epsilon(x) = \epsilon(\beta 1 + \gamma x) = \beta \epsilon(1) + \gamma \epsilon(x) = \gamma a.
\]
Multiplying by $f$ and recalling that $af=1$ we see $\gamma=1$. 

By axiom (1) of \fullref{def:efa} we have $\theta = \phi(\theta)$ and $\theta x=\phi (\theta x)$. These two relations imply that $\mu \beta = 0$ and $\lambda \beta =0$.

By axiom (2) we have $m(\phi\otimes \Id)((\Delta(1)) = \theta^2$. Now the left-hand side is equal to $(f\beta - hf) 1 + 2f x$ and the right-hand side is equal to $ (\lambda^2 + \mu^2 t)1 + (2\lambda\mu + \mu^2h)x$. 
Thus we have
\[
\lambda^2 + \mu^2 t = f\beta - hf \quad \text{and} \quad 2\lambda\mu + \mu^2h = 2f.
\]
Recalling that $af=1$ the first equation gives $h=\beta - a\lambda^2
- a\mu^2t$ and the second $2a\lambda\mu + a \mu^2h = 2$. Substituting for $h$ in the latter gives
\[
2a\lambda\mu - a^2\mu^2\lambda^2 -a^2\mu^4t = 2.
\]
The additional demand that $\phi$ is an involution and respects
multiplication and comultiplication does not introduce any further
relations.

Thus $\theta$ must be of the form
$ \theta = \lambda 1 + \mu x $ and $\phi$ must be of the form
$ \phi(1) = 1$ and  $\phi(x) = \beta 1 + x$ where
$\lambda, \mu, \beta\in U$ satisfy
\[
\mu\beta= \lambda\beta = 0 \quad \text{and} \quad 2a\lambda\mu - a^2\mu^2\lambda^2 -a^2\mu^4t = 2,
\]
and there are no further relations necessary. Along with the relation
$af=1$ these are the relations defining $K$ so there is a map
$\psi\colon K\ra R$ as required. Noting that $ h= \beta - a \lambda^2
-a \mu^2t$ one easily sees that $U_\psi \cong V$.
\end{proof}


\section{The complex of cobordisms for a link diagram on a surface}\label{sec:complex}

\begin{figure}[ht!]
\labellist
\small\hair 2pt
\pinlabel $1$ [br] at 45 85
\pinlabel $2$ [b] at 124 75
\pinlabel $3$ [bl] at 170 79
\endlabellist
\centerline{\includegraphics{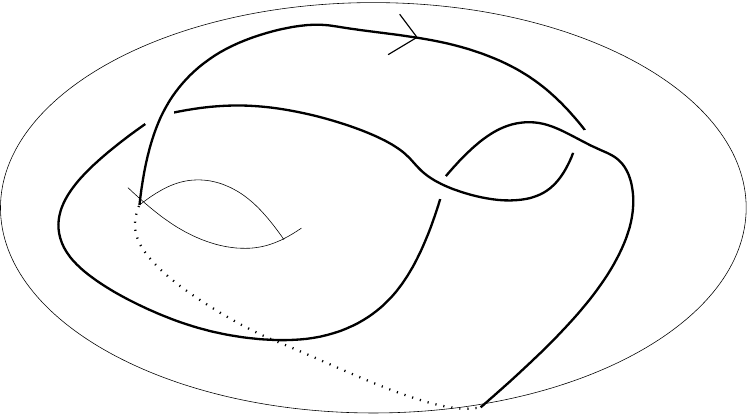}}
\caption{An oriented knot diagram on a torus}
\label{fig:32knot} 
\end{figure}

We consider oriented link diagrams on closed orientable surfaces (as
depicted for example in \fullref{fig:32knot}) and define an
equivalence relation on such diagrams as follows. Two diagrams are
said to be {\em stably equivalent\/} if they are related by
\begin{enumerate}
\item surface homeomorphisms
\item Reidemeister moves
\item addition or subtraction of handles to the surface (when 
  no part of the diagram is on the handle).
\end{enumerate}

Given a link diagram on a closed orientable surface $\Sigma$ we wish to
define a complex of cobordisms along the lines of Bar-Natan's complex
for link diagrams given in \cite{barnatan}. 

Recall that for an
additive category $\cC$ the additive category $\mat(\cC)$ is defined as
follows. Its objects are finite families $\{C_i\}_i$ of objects
$C_i\in \cC$ which for convenience will be denoted $\oplus C_i$. A
morphism $F\colon \oplus C_i \ra \oplus C^\prime_{\smash j}$ is a matrix
$F=[F^{\smash[t]{j}}_{\smash[b]{i}}]$ of morphisms $F^{\smash[t]{j}}_{\smash[b]{i}}\colon C_i \ra C^\prime_{\smash j} $ in
$\cC$. We will refer to the $F^{\smash[t]{j}}_{\smash[b]{i}}$ as the {\em matrix elements\/}
of the morphism $F$. Composition is defined in terms of matrix
elements by the rule $[F \circ G]^{\smash k}_{\smash[b]{i}} = \sum_j F^{\smash k}_{j}\circ
G^{\smash[t]{j}}_{\smash[b]{i}}$ and addition of morphisms given by matrix addition.  If
$\cC$ is not additive then it is made so by allowing formal
$\bZ$--linear combinations of morphisms.

Define $\ucob \Sigma$ to be the following category. The objects are
collections of disjoint closed curves $\Gamma$ in $\Sigma$. A morphism
$\Gamma\ra \Gamma^\prime$ is a surface embedded in $\Sigma \times I$
whose boundary lies entirely in $\Sigma\times \{0,1\}$ and which
agrees with $\Gamma$ on $\Sigma\times \{0\}$ and with $\Gamma^\prime$
on $\Sigma\times \{1\}$. Two such morphisms are identified if they are
related by a boundary preserving isotopy. Since cobordisms are
embedded one can clearly compose them (rescaling the result).

\begin{figure}[ht!]
\begin{center}
\includegraphics{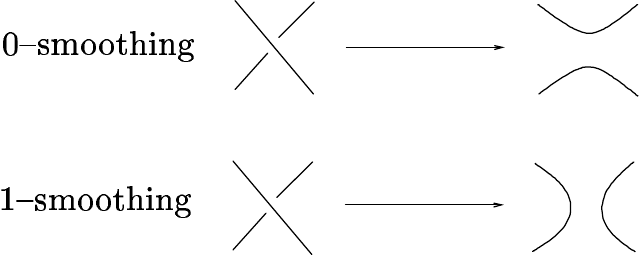}
\end{center}
\caption{}
\label{fig:smoothings}
\end{figure}

Given an oriented link diagram $D$ on an orientable surface $\Sigma$ we first
number the crossings $1,\ldots, n$. We can resolve each crossing in
one of two ways as depicted in \fullref{fig:smoothings}.

\begin{figure}[ht!]
\labellist
\tiny\hair 0pt
\pinlabel $100$ [t] at 117 138
\pinlabel $101$ [t] at 204 138
\pinlabel $000$ [t] at 22 72
\pinlabel $010$ [t] at 118 71
\pinlabel $101$ [t] at 204 73
\pinlabel $111$ [t] at 288 73
\pinlabel $001$ [t] at 118 4
\pinlabel $011$ [t] at 203 6
\endlabellist
\centerline{\includegraphics{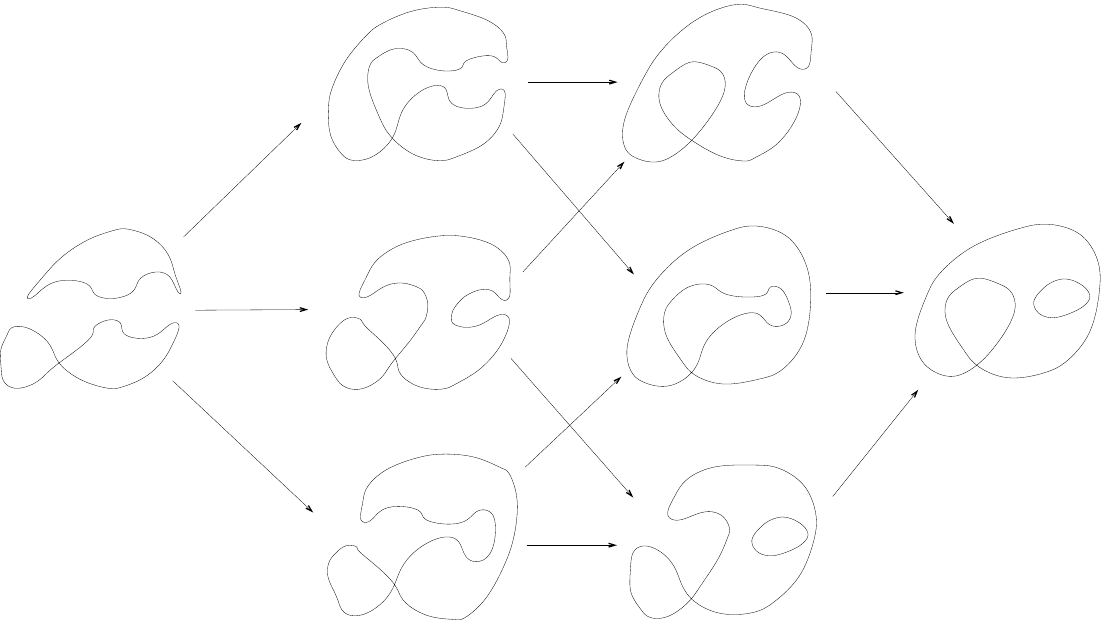}}
\caption{The cube of the diagram in \fullref{fig:32knot}}
\label{fig:poset32} 
\end{figure}

A resolution of each of the $n$ crossings of $D$ will be called a {\em
smoothing\/}. There are $2^n$ such smoothings and each is indexed by a
sequence $s$ of $n$ 0's and 1's, the $i$--th entry informing us whether
the $i$--th crossing is a $0$-- or $1$--smoothing. The smoothing itself
consists of a closed $1$--manifold $\Gamma_s$ being a collection of
nonintersecting closed curves in $\Sigma$. Note that $\Gamma_s$ is an
object in the category $\ucob \Sigma$. Smoothings form a poset via
$$s<t \text{ iff all $1$--smoothings in $s$ are $1$--smoothings in $t$}.$$  
\begin{gather*}
r(s) = \text{number of $1$--smoothings in $s$} \tag*{\hbox{Let}}\\
k(s) = \text{number of components in $\Gamma_s$}.\tag*{\hbox{and}}
\end{gather*}
As is familiar the smoothings are arranged on the vertices of a cube with
an arrow from $s$ to $t$ if $s<t$ and $r(t) = r(s) +1$. The cube for
the knot given in \fullref{fig:32knot} is presented in \fullref{fig:poset32} where the underlying torus has been omitted.

If $s$ and $t$ are smoothings such that $r(t)= r(s)+1$ and $s<t$ then
$\Gamma_t$ must be identical to $\Gamma_s$ outside a small disc in
$\Sigma$ around a $1$--smoothing in $\Gamma_t$. We refer to this disc as the
{\em changing disc\/}. In this situation we define $\langle s,t\rangle $
to be the number of $1$--smoothings among the first $j-1$ crossings of
$t$ where the $j$--th crossing is the one in the changing disc.

For such $s$ and $t$ define $W_s^t\subset \Sigma\times I$ to be the
surface which is $\Gamma_s$ on $\Sigma\times \{0\}$, $\Gamma_1$ on
$\Sigma\times
\{1\}$, a product outside $(\text{changing disc})\times I$ and a
saddle in place of the missing $(\text{changing disc})\times I$. We regard $W_s^t$ as a morphism $\Gamma_s \ra \Gamma_t$ in the category $\ucob \Sigma$. 
\[
\bri D i = \bigoplus_s \Gamma_s \in \text{Ob}(\mat(\ucob \Sigma)) \leqno{\hbox{Set}}
\]
where the sum is over all smoothings $s$ with $r(s) = i+n_-$. Here
$n_-$ is the number of negative crossings in $D$.

We now want to define a morphism $d^i\colon \bri D i \ra \bri D {i+1}$.
In order to define
$d^i$ it is enough to define its matrix elements $(d^i)^{t}_{s}\colon
\Gamma_s\ra \Gamma_t$ where $s$ and $t$ are smoothings such that
$r(s) = i+n_-$ and $r(t)= i+1+n_-$. Define these matrix elements by
\[
(d^i)^{t}_{s} = \begin{cases}
(-1)^{\langle s,t\rangle} W^{t}_{s} & \text{ if $s<t$}\\
0 & \text{else.}
        \end{cases}
\]
\begin{prop}
Given a diagram $D$, the morphism $d$ defined  above satisfies
$d^2=0$. Thus $(\bri D *, d)$ is a complex in $\mat(\ucob \Sigma)$.
\end{prop}
\begin{proof}
Given a square face of the cube
\[
\xymatrix{
& t \ar[dr] &  \\
s \ar[ur] \ar[dr] & & u \\
& t^\prime \ar[ur]  & 
}
\]
we note that $W_t^u \circ W_s^t\cong W_{t^\prime}^u
\circ W_s^{t^\prime}$ since saddles can be reordered. 
The signs chosen in the definition of the matrix elements ensure that
each square anticommutes.
\end{proof}

\begin{figure}[ht!]
\labellist
\small\hair 2pt
\pinlabel $=0$ [l] at 26 32
\pinlabel ${\it S\/}$ [t] at 26 8
\pinlabel $=2$ [l] at 92 32
\pinlabel ${\it T\/}$ [t] at 89 8
\pinlabel* $+$ at 164 31
\pinlabel* $=$ at 206 31
\pinlabel* $+$ at 245 31
\pinlabel $\textit{4--Tu}$ [t] at 205 8
\endlabellist
\centerline{\includegraphics{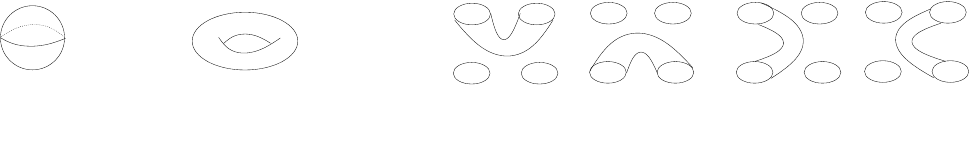}}
\caption{Bar-Natan's relations}
\label{fig:barnatanrels} 
\end{figure}
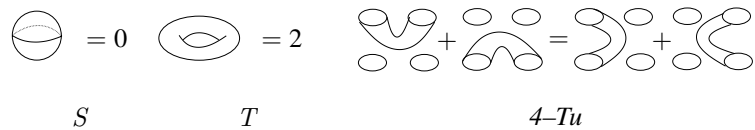

We continue following \cite{barnatan} by quotienting the category
$\ucob \Sigma$ by certain relations. These relations, the {\em S,T,
4--Tu\/} (Sphere, Torus and $4$--Tube) relations are illustrated below in
\fullref{fig:barnatanrels}.

\begin{defn}
Let $\ucobmodr \Sigma$ be the category obtained from $\ucob \Sigma$ by
quotienting by the equivalence relation generated by relations {\em
S,T,4-Tu\/}.
\end{defn}

\begin{prop}\label{prop:invariance} 
Let $D$ be a link diagram on a closed orientable surface $\Sigma$. The
homotopy type of the complex $(\br D, d)$ in $\mat (\ucobmodr \Sigma)$ is invariant under Reidemeister moves.
\end{prop}

\begin{proof}
 In \cite{barnatan} Bar-Natan has proved the analogous statement for
 planar link diagrams. His proofs remain
 valid in our setting too.
\end{proof}


\section{Link homology}\label{sec:homology}
\subsection{Applying an unoriented TQFT to get link homology}
In this section we will apply a rank-two aspherical unoriented TQFT to
the formal complex of cobordisms of the previous section to obtain a
complex of $R$--modules. We can then take homology of this complex to
obtain a calculable invariant of stable equivalence classes of link
diagrams on surfaces.

While geometrically the {\em 4-Tu\/} relation is a genuine relation on cobordisms, algebrai\-cally it comes for free for aspherical theories.

\begin{prop}\label{prop:4tu}
Let $(A,\tau)$ be a rank-two aspherical unoriented TQFT. Suppose
cobordisms $W_1,W_2,W_3,W_4$ are related locally as in the {\it 4-Tu\/}
relation above then $\tau(W_1) + \tau(W_2) = \tau(W_3) + \tau(W_4)$.
\end{prop}

\begin{proof}
In order to do computations we need to work with the underlying
Frobenius algebra and in order to do that we need to pick
identifications $\gamma_i\colon V\ra \Gamma_i$, $i=1,2,3,4$, where the 
$\Gamma_i$ are the four circles appearing $W_1$ (and $W_2,W_3,W_4$) where the local surgery takes place. (Note that we are
assuming that $W_1,W_2,W_3,W_4$ are actually the same manifold outside
the ball in which the local change takes place.)
We need to show:

\begin{figure}[ht!]
\centerline{\includegraphics{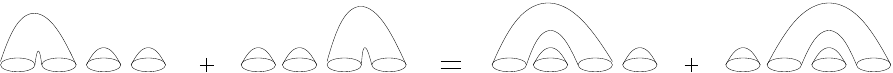}}
\end{figure}
which by using the identifications $\gamma_i$ can be converted into
the condition on the underlying Frobenius algebra:
\[
\Sigma a^\prime \ot a^{\prime\prime}\ot 1\ot 1+
\Sigma 1\ot 1\ot a^\prime \ot a^{\prime\prime} =
\Sigma a^\prime \ot 1\ot a^{\prime\prime}\ot 1+
\Sigma 1\ot a^\prime \ot 1 \ot a^{\prime\prime}
\]
where we write $\Delta(1) =\Sigma a^\prime \ot a^{\prime\prime}$.
By \fullref{prop:univ} it suffices to consider the
universal extended Frobenius algebra. Since here
$\Delta(1) = f(1\ot x +x\ot 1)- hf 1\ot 1$ it is easy to verify that
the above equation holds.
\end{proof}

Let $(A,\tau)$ be a rank-two aspherical unoriented TQFT. We can now
define a functor $\ucob \Sigma \ra \text{Mod}_R$ by
\[
\Gamma \mapsto A_{\Gamma},
\]
\[
W \mapsto \tau(W).
\]

This does indeed give a functor since by axiom (4) in \fullref{def:tqft} homeomorphisms preserving the boundary point-wise induce the
same map $\tau$ and composition in $\ucob \Sigma$ glues along the
identity map and so $\tau(W^\prime \circ W) = \tau(W^\prime) \circ \Id
\circ
\tau(W)$. Note also that by axiom (7) cylinders give the identity
homomorphism.

By taking formal direct sums to genuine direct
sums this extends to a functor
\[
\mat(\ucob \Sigma) \ra \modR.
\]
Since $(A,\tau)$ is rank-two and aspherical, \fullref{prop:4tu} shows that this functor factors through $\mat(\ucobmodr \Sigma)$ and hence there is 
a functor $F_{(A,\tau)}$ on the associated categories of complexes
\[
F_{(A,\tau)}\colon \kom(\mat(\ucobmodr \Sigma)) \ra \kom(\modR).
\]
Given a link diagram $D$ on a surface  set
\[
C^{*}_{(A,\tau)}(D) = F_{(A,\tau)}(\bri D*).
\]
\begin{prop}
The homotopy type of the complex $C^{*}_{(A,\tau)}(D)$ is an invariant of stable equivalence classes.
\end{prop}

\begin{proof}
Let $\Sigma\ra \Sigma^\prime$ be a homeomorphism of surfaces. This
induces homeomorphisms $\Gamma_s \ra \Gamma^\prime_s$ on smoothings
and homeomorphisms $W_s^t \ra (W^\prime)_s^t$ on the cobordisms defining the
complexes. Applying the TQFT gives isomorphisms $A_{\Gamma_s}
\ra A_{\Gamma^\prime_s}$ which in turn induce an isomorphism (of graded vector spaces)
$C^{*}_{(A,\tau)}(D) \ra
C^{*}_{(A,\tau)}(D)^\prime$. By axiom (4) in the definition
of TQFT we see that this is an isomorphism of complexes.

Showing invariance under Reidemeister moves is the content of \fullref{prop:invariance}.

If $\Sigma^\prime$ is obtained from $\Sigma$ by the addition or
subtraction of handles disjoint from the diagram then clearly there is a
canonical identification of smoothings and cobordisms in the complexes
${\br D}^\prime$ and $\br D$ (defined above using $\Sigma^\prime$ and
$\Sigma$). These identifications lead to an isomorphism of complexes $C^{*}_{(A,\tau)}(D) \ra
C^{*}_{(A,\tau)}(D)^\prime$.
\end{proof}

Thus if $L$ is a stable equivalence class of diagrams we can now define what we mean by link homology.
\begin{defn}
The {\it link homology for stable equivalence classes based on the TQFT $(A,\tau)$\/} is defined by 
\[
H^*_{(A,\tau)}(L) = H(C^*_{(A,\tau)}(D))
\]
where $D$ is any representative diagram of the class $L$.
\end{defn}

Note that if we replace the unoriented TQFT $(A,\tau)$ with an
isomorphic one $(A^\prime,\tau^\prime)$ then there are
isomorphisms $\eta\colon A_{\Gamma_s} \ra A^\prime_{\Gamma_s}$. As
these are natural with respect to cobordisms they induce an
isomorphism of complexes $C^{*}_{(A,\tau)}(D) \ra
C^{*}_{(A^\prime,\tau^\prime)}(D)$. Thus isomorphic TQFTs result in
isomorphic link homology groups.

We now present our main result.

\begin{thm}\label{thm:vlh}
 Let $R$ be a commutative ring with unit. Let  $\lambda,
\mu,\beta, a, t\in R$ such that $a$ is invertible and the following relations are satisfied:
\begin{eqnarray}
 &\mu\beta= \lambda\beta = 0,&\label{eq:1}\\ 
& 2a\lambda\mu - a^2\mu^2\lambda^2
 -a^2\mu^4t = 2.&\label{eq:2}
\end{eqnarray}
Then the $5$--tuple $(\lambda, \mu,\beta,a, t)$ defines (the isomorphism
class of) a link homology theory for stable equivalence classes of
diagrams on oriented surfaces.
\end{thm}

\begin{proof}
The obvious map $\psi\colon K\ra R$ defines a rank-two extended
Frobenius algebra $(V,\theta,\phi)$ satisfying $\epsilon(i(1))
=0$. This in turn defines a rank-two aspherical $1{+}1$--dimensional
unoriented TQFT. Such a theory defines a link homology for stable
equivalence classes as described above.
\end{proof}

\begin{cor}
If $R$ is an integral domain then each triple $(a,\lambda, \mu)$ with  $a$ and $\mu$ invertible
defines a link homology.
\end{cor}

\subsection[Link homologies over Q]{Link homologies over $\bQ$}
Clearly there are solutions to equations \eqref{eq:1} and \eqref{eq:2}
over $\bQ$. The question we now address is whether or not any of these
extend a theory isomorphic to Khovanov's original link homology. Given
$h,t\in \bQ$ then there is a link homology defined by the Frobenius
algebra $\bQ\{1,x\}$ with $i(1)=1,\epsilon(1) =0,\epsilon(x) =1$ and
$x^2=hx+t1$. Recall from \cite{mackaayturner} that such a theory is
isomorphic to Khovanov's original theory if and only if
$h^2+4t=0$. Noting that over $\bQ$ we must have $\beta=0$ and
$\mu\neq 0$ we can express $h$ and $t$ in terms of $\lambda$ and $\mu$
as
\[
h=2\mu^{-2} - 2\lambda\mu^{-1}
\quad \text{and} \quad
t=-\mu^{-2}\lambda^2 -2\mu^{-4} + 2\lambda\mu^{-3}.
\]
A computation now shows
\[
h^2+4t=-4\mu^{-4}.
\]
Thus we have $h^2+4t \neq 0$ resulting in the disappointing conclusion
that Khovanov's original theory does not extend (at least using the
methods of this paper) to stable equivalence classes. This conclusion remains valid over fields of characteristic not equal to two.

Over $\bQ$ the resulting theories are all singly
graded. For each of these the Euler characteristic, $\chi$, is the unnormalised
Jones polynomial evaluated at $q=1$. To see this recall that the Euler
characteristic of the homology of a complex is the same as the Euler
characteristic of the complex itself and so
\begin{align*}
\chi = \sum_i (-1)^i & \text{dim}(H^i(D)) =  \sum_i (-1)^i \text{dim}(F_{(A,\tau)}(\bri D i)) \\
& = \sum_i (-1)^i \sum_{\substack{s\in \text{\small smoothings} \\
r(s) = i+n_-}}\text{dim}( A_{\Gamma_s}) = \sum (-1)^i
\sum_{\substack{s\in \text{\small smoothings} \\ r(s) = i+n_-}}
2^{k(s)}
\end{align*}
where $k(s)$ is the number of components in $\Gamma_s$.
Comparing with the definition of the Jones polynomial we see that
the right-hand side is the unnormalised Jones polynomial evaluated at 1.

\subsection[Link homologies over F_2]{Link homologies over $\FF$}
Working over $\FF$ one has more success than over the rationals. One has $a=1$ and the remaining equations of the main theorem become
$$
\mu\beta= \lambda\beta = 0
\quad \text{and} \quad
\mu^2\lambda^2 +\mu^4t = 0.
$$
There are eight possibilities  as tabulated below.
\begin{center}\small
\begin{tabular}{|c|c|c|c|c|c|c|}
\hline
\strut$\lambda$ & $\mu$ & $t$ & $\beta$ &$h$&$\theta$&$\phi(x)$\\\hline
\strut0&0&0&0&0&0&$x$\\
\strut0&0&0&1&1&0&$1+x$\\
\strut1&0&0&0&1&1&$x$\\
\strut0&0&1&0&0&0&$x$\\
\strut0&0&1&1&1&0&$1+x$\\
\strut1&0&1&0&1&1&$x$\\
\strut0&1&0&0&0&$x$&$x$\\
\strut1&1&1&0&0&$1+x$&$x$\\\hline
\end{tabular}
\end{center}

Rows 1, 4, 7 and 8 have isomorphic underlying Frobenius algebras and on
classical links give theories that are isomorphic to Khovanov's
original theory over $\FF$. As extended Frobenius algebras there are
two isomorphism classes among these with row 1 isomorphic to row 4 and
row 7 isomorphic to row 8.  The theory in the first row gives a
bigraded theory and was investigated by Manturov \cite{manturov}. He
observes that the graded Euler characteristic (which can be computed
from the graded Euler characteristic of the chain complex) is the
unnormalised Jones polynomial of the link.

The theory in row 3 can be made into a bigraded theory by taking  $R=\FF[\lambda]$ with $\lambda$ in degree
$-1$. The Frobenius algebra  has
multiplication
\[
11=1, \quad 1x=x1=x \quad \text{and} \quad xx= \lambda^{2}x 
\]
and comultiplication
\begin{gather*}
\Delta (1) = 1\otimes x + x\otimes 1 + \lambda^{2} 1\otimes 1,\\
\Delta(x) = x\otimes x. 
\end{gather*}
The counit is $\epsilon (1) = 0,\epsilon (x) =1$. This is 
isomorphic to Bar-Natan's graded characteristic two theory. For the
extended structure we have $\theta = \lambda 1$and $\phi=\Id$. This theory is
related to Manturov's theory: one can filter the chain complex by
powers of $\lambda$ to produce a spectral sequence (similar to that in
\cite{turner}) whose $E_2$--page is Manturov's theory and which
converges to the bigraded theory.

It is interesting to note that this extended Frobenius algebra has an interpretation in terms of the equivariant cohomology of $S^2$ along the lines of \cite{khovanov_frob}.
Let $G=\bZ/2$ act
on $S^2$ by a rotation through $\pi$ about a fixed axis. In this case
we have
\[
H^*_G(\pt;\FF) = H^*(\bR P^\infty;\FF) = \FF[\lambda]
\] 
\[
H^*_G(S^2;\FF) = \FF[\lambda,x]/\langle x^2 = \lambda^2x\rangle = \FF[\lambda][x]/\langle x^2 = \lambda^2x\rangle .
\leqno{\hbox{and}}
\]
We see that this is the  algebra of the bigraded link
homology above and the element $\theta$ is the image of the generator
of $H^*_G(\pt;\FF)$ under the homomorphism induced from $S^2 \ra \pt$.

Bigraded Bar-Natan theory is also extended by the theory in row
2. Here we work over $R=\FF[\beta]$ with $\beta$ in degree $-2$. In
this case $\theta=0$ and $\phi(x) = \beta 1+x$. Note that in this case we
have a nontrivial involution $\phi$.

\subsection{Acknowledgments} 
We thank D\,Bar-Natan and G\,Naot for many helpful comments. The second
author was supported by the European Commission through a Marie Curie
fellowship and thanks the Institut de Recherche Math\'ematique
Avanc\'ee in Strasbourg for their hospitality.


\bibliographystyle{gtart}
\bibliography{link}

\end{document}